\documentclass[12pt]{article}
\usepackage{graphicx}
\usepackage{amsmath}
\usepackage{amsfonts}
\usepackage{amssymb}
\newtheorem{theorem}{Theorem}[section]
\newtheorem{proposition}[theorem]{Proposition}
\newtheorem{lemma}[theorem]{Lemma}
\newtheorem{corollary}[theorem]{Corollary}

\begin{document}

\title{Relaxed Gabor expansion at critical density and a `certainty principle'}
\author{V.P.Palamodov\\Tel Aviv University}
\date{31.07.05}
\maketitle

\section{Introduction}

The functions in $\mathbb{R}^{n}$ that are the best localized in configuration
and in frequency spaces simultaneously are known in quantum field theory as
coherent states (shifted modulated gaussian functions). Any coherent state
gives exact minimum for the uncertainty relation which means that each
coherent state $\mathbf{e}_{\lambda}$ is maximally concentrated at a certain
point $\lambda$ of the phase space $\Phi=\mathbb{R}^{n}\times\left(
\mathbb{R}^{n}\right)  ^{\prime}$. Another advantage of the family $\left\{
\mathbf{e}_{\lambda},\lambda\in\Phi\right\}  $ is its large symmetry group:
the Weyl-Heisenberg group is acting by shifts and modulation operators and the
metaplectic group which is a `quantization' of the group of linear simplectic transformations.

The role of the Weyl-Heisenberg family in the information theory
($n=1$) was emphasized by Gabor \cite{Ga}: ''each elementary
signal conveys exactly one datum, or one quantum of information''.
Gabor's idea was to expand an `arbitrary' signal in a series of
coherent states $\mathbf{e}_{\lambda}$ for points $\lambda$ in a
maximally sparse lattice $\Lambda$ in the phase plane. The
property of the Gabor system $\left\{ \mathbf{e}_{\lambda},\lambda
\in\Lambda\right\}  $ depends on the area $a$ of a cell of the
lattice $\Lambda.$ The system is complete, if and only if
$a\leq1$. Lyubarskii \cite{L} has shown that the Gabor system is a
frame, if $a<1$. We focus on the critical case $a=1$ when the
Gabor system of critical density is still
complete, but is not a frame. A function $f\in$ $L_{2}\left(  \mathbb{R}%
\right)  $ may not have a convergent Gabor series, but there
exists always a series that converges in a distribution sense, see
Janssen \cite{Ja}. The book of I.Daubechies \cite{Da} contains a
survey of the problem, see also Feichtinger, Strohmer \cite{FS}
for applications of the Gabor analysis.

Lyubarskii and K.Seip \cite{LS} has shown that for an arbitrary
function $f\in L_{2}$ that is 1/2-smooth in time and frequency
there exists a convergent series over a Gabor system, if a minor
deformation of a lattice $\Lambda$ of critical density is made.

We show here that any $\delta$-smooth (in $\Phi$) function for
$\delta>1$ has a \textit{unique} convergent Gabor expansion with
$l_2$-coefficients, provided one more point $\sharp$ is added to a
lattice $\Lambda $; this point is in the middle of a cell and is
called \textit{sharp} point. We call it a \textit{relaxed} Gabor
expansion. This approach can save Gabor's idea. The rate of
convergence of the coefficients can be improved, if few more sharp
terms are added to the Gabor system. The coefficients of such a
series are again uniquely defined for a sufficiently smooth
function $f$.

The relaxed Gabor expansion sheds light to the physical wisdom:
\textit{a field supported in certain domain} $D$ \textit{in the
phase space has about }$\left|  D\right|  $ \textit{degrees of
freedom}, where $\left|  D\right|  $ means the simplectic area of
the domain (Nyquist, Wigner, Brillouin, Shannon, Gabor,...). This
claim can be called `certainty' principle as opposite of the
classical uncertainty one. The Landau-Pollak dimension theorem
\cite{LP} gives an accurate form to this principle in a special
situation. Given positive numbers $T,\Omega,$ any function $f$
`concentrated' in the box $D$ of the size $T\times2\Omega$ can be
approximated by a linear combination of $2T\Omega +O\left(  \log
T\Omega\right)  $ first prolate spheroidal functions. The
condition of concentration means that the tail of a unit energy
function $f$ is small in $\mathbb{R\backslash}\left[  0,T\right] $
and the tail of Fourier transform $\hat{f}$ is small in
$\mathbb{R}^{\ast}\backslash\left[ -\Omega,\Omega\right]  .$ The
number $2T\Omega$ is equal to the area $\left| D\right|  $ of the
box measured by the canonical simplectic form in the phase plane.
This means the approximation `dimension' of the space of functions
localized in the box equals the simplectic area of this up to a
logarithmic term.

We formulate and prove here a rigorous form of the certainty principle for
domains of arbitrary shape. We use Gabor means $\left\langle f|\mathbf{e}%
_{\mu}\right\rangle $ to specify the condition of concentration of a function
$f\in L_{2}$ in the phase space $\Phi$. We say that $f$ is concentrated in a
domain $D$, if the integral over $\Phi\backslash D$ of the density $\left|
\left\langle f|\mathbf{e}_{\mu}\right\rangle \right|  ^{2}\mathrm{d}\mu$ is
small. We show that for an arbitrary sets $K\subset D,$ where $D$ is
$r$-neighborhood of $K$ for some $r\geq r_{0}$, any function $f$ concentrated
in $D$ can be written in the form $g+h+\phi,$ where $g$ is a linear
combination of Gabor functions $\mathbf{e}_{\lambda}$, with $\lambda\in
\Lambda\cap D$, $h$ belongs to the linear envelope of Gabor functions with
sharp $\lambda\in D\backslash K$ and the norm of $\phi$ is bounded by the rate
of concentration of $f$ plus an exponentially small term as $r\rightarrow
\infty$. The total number $N$ of such points $\lambda$ is equal to $\left|
D\right|  +\left|  D\backslash K\right|  =\left|  D\right|  +O\left(  r\left|
D\right|  ^{1/2}\right)  $. To summarize briefly, this means that any function
$f$ concentrated in $D$ is a linear combination of $N\sim\left|  D\right|
+O\left(  r\left|  D\right|  ^{1/2}\right)  $ coherent states supported in $D$
up to a exponentially small term.

I thank Yu.Lyubarskii for helpful discussions.

\section{Gabor transform and localization}

Fix a coordinate $x$ in a line $\mathbb{R}$; the family of functions in
$\mathbb{R}$
\[
\mathbf{e}_{\lambda}\left(  x\right)  \doteq2^{1/4}\exp\left(  -\pi\left(
x-p\right)  ^{2}+2\pi\imath\theta x\right)
\]
are called Gabor functions. Here $\lambda=\left(  p,\theta\right)  $ is a
point in the phase space $\Phi\doteq\mathbb{R}\times\mathbb{R}^{\prime}$
($\mathbb{R}^{\prime}$ means the dual line). Another term for $\mathbf{e}%
_{\lambda}$ is \textit{coherent state} at the point $\lambda$. The phase space
has the natural Euclidean structure:\ $\left|  \lambda\right|  ^{2}=\left|
p\right|  ^{2}+\left|  \theta\right|  ^{2}$ and Lebesgue measure
\textrm{d}$\lambda\doteq\mathrm{d}p\,\mathrm{d}\theta.$ We use the notation
$\left\langle \cdot|\cdot\right\rangle $ for the scalar product in
$L_{2}=L_{2}\left(  \mathbb{R}\right)  $ and $\left\|  \cdot\right\|  $ for
the norm. For $f\in L_{2}$ the function $\lambda\mapsto\left\langle
f|\mathbf{e}_{\lambda}\right\rangle $ defined in $\Phi$ is called
\textit{Gabor transform} of $f.$ This transform is unitary:

\begin{proposition}
\label{fee}For an arbitrary function $f\in L_{2}$ the equation holds
\begin{equation}
\int_{\Phi}\left|  \left\langle f|\mathbf{e}_{\lambda}\right\rangle \right|
^{2}\mathrm{d}\lambda=\left\|  f\right\|  ^{2}\label{fp}%
\end{equation}
where $\mathrm{d}\lambda=\mathrm{d}p\mathrm{d}\theta$ and
\begin{equation}
f=\int_{\Phi}\left\langle f|\mathbf{e}_{\lambda}\right\rangle \mathbf{e}%
_{\lambda}\mathrm{d}\lambda\label{fed}%
\end{equation}
in weak $L_{2}$-sense.
\end{proposition}

\textsc{Proof.} By Plancherel Theorem
\begin{align}
\int_{\mathbb{R}}\left|  \left\langle f|\mathbf{e}_{\lambda}\right\rangle
\right|  ^{2}\mathrm{d}\theta &  =2^{1/2}\int\left|  \int f\left(  x\right)
\exp\left(  -\pi\left(  x-p\right)  ^{2}\right)  \exp\left(  2\pi\imath\theta
x\right)  \mathrm{d}x\right|  ^{2}\mathrm{d}\theta\nonumber\\
&  =2^{1/2}\int\left|  F\left(  f_{p}\right)  \left(  -\theta\right)  \right|
^{2}\mathrm{d}\theta=2^{1/2}\int\left|  f_{p}\left(  x\right)  \right|
^{2}\mathrm{d}x, \label{dt}%
\end{align}
where $f_{p}\left(  x\right)  =f\left(  x\right)  \exp\left(  -\pi\left(
x-p\right)  ^{2}\right)  .$ Next we integrate both sides against $\mathrm{d}p$
and changing the variable $x$ to $y=x-p.$ This yields
\begin{align*}
\int\int\left|  f_{p}\left(  x\right)  \right|  ^{2}\mathrm{d}x\mathrm{d}p  &
=\int\int\left|  f\left(  x\right)  \exp\left(  -\pi y^{2}\right)  \right|
^{2}\mathrm{d}x\mathrm{d}y\\
&  =\int\left|  f\left(  x\right)  \right|  ^{2}\mathrm{d}x\int\exp\left(
-2\pi y^{2}\right)  \mathrm{d}y\\
&  =2^{-1/2}\int\left|  f\left(  x\right)  \right|  ^{2}\mathrm{d}x
\end{align*}
and (\ref{fp}) follows. Now calculate the scalar product of integral
(\ref{fed}) with $\mathbf{e}_{\mu}.$ We have for arbitrary points
$\lambda=\left(  p,\theta\right)  ,\mu=\left(  q,\eta\right)  \in\Phi$%
\begin{equation}
\left\langle \mathbf{e}_{\lambda}|\mathbf{e}_{\mu}\right\rangle =\exp\left(
\pi\imath\left(  p+q\right)  \left(  \theta-\eta\right)  -\pi\left|
\lambda-\mu\right|  ^{2}/2\right)  . \label{ee}%
\end{equation}
Therefore%
\[
\int_{\Phi}\left\langle f|\mathbf{e}_{\lambda}\right\rangle \left\langle
\mathbf{e}_{\lambda}|\mathbf{e}_{\mu}\right\rangle \mathrm{d}\lambda
=\left\langle f|E_{\lambda,\mu}\right\rangle ,
\]
where
\begin{align*}
\bar{E}_{\lambda,\mu}  &  \doteq\int_{\Phi}\left\langle \cdot|\mathbf{e}%
_{\lambda}\right\rangle \left\langle \mathbf{e}_{\lambda}|\mathbf{e}_{\mu
}\right\rangle \mathrm{d}\lambda\\
&  =\int_{\Phi}\left\langle \cdot|\mathbf{e}_{\lambda}\right\rangle
\exp\left(  \pi\imath\left(  p+q\right)  \left(  \theta-\eta\right)
-\pi\left|  \lambda-\mu\right|  ^{2}/2\right)  \mathrm{d}\lambda\\
&  =2^{1/4}\int\int\exp\left(  -\pi\left(  x-p\right)  ^{2}-2\pi\imath\theta
x\right)  \exp\left(  \pi\imath\left(  p+q\right)  \left(  \theta-\eta\right)
-\pi\left|  \lambda-\mu\right|  ^{2}/2\right)  \mathrm{d}\theta\mathrm{d}p
\end{align*}
Calculate the interior integral by changing $\theta$ to $\xi=\theta-\eta$:%
\begin{align*}
\bar{E}_{\lambda,\mu}\left(  x\right)   &  =2^{1/2}\exp\left(  -2\pi\imath\eta
x\right)  \int\exp\left(  -\pi\left(  x-p\right)  ^{2}-\frac{\pi}{2}\left(
p+q-2x\right)  ^{2}-\frac{\pi}{2}\left|  p-q\right|  ^{2}\right)
\mathrm{d}p\\
&  =2^{1/2}\exp\left(  -2\pi\imath\eta x\right)  \exp\left(  -\pi\left(
x-q\right)  ^{2}\right)  \int\exp\left(  -2\pi\left(  x-p\right)  ^{2}\right)
\mathrm{d}p=\mathbf{\bar{e}}_{\mu}\left(  x\right)
\end{align*}
This allows to write%
\[
\int_{\Phi}\left\langle f|\mathbf{e}_{\lambda}\right\rangle \left\langle
\mathbf{e}_{\lambda}|\mathbf{e}_{\mu}\right\rangle \mathrm{d}\lambda
=\left\langle f|\mathbf{e}_{\mu}\right\rangle
\]
This yields for an arbitrary $g\in L_{2}$%
\[
\left\langle f|g\right\rangle =\lim_{r\rightarrow\infty}\int_{\left|
\lambda\right|  \leq r}\left\langle f|\mathbf{e}_{\lambda}\right\rangle
\left\langle \mathbf{e}_{\lambda}|g\right\rangle \mathrm{d}\lambda=\int_{\Phi
}\left\langle f|\mathbf{e}_{\lambda}\right\rangle \left\langle \mathbf{e}%
_{\lambda}|g\right\rangle \mathrm{d}\lambda
\]
The integral converges, since both factors belong to $L_{2}\left(
\Phi\right)  .$ $\blacktriangleright$

Now we show that a function $f\in L_{2}$ can be localized in a convex domain
$D\subset\Phi$ in terms of its Gabor transform.

\begin{corollary}
\label{hp}For an arbitrary $f\in L_{2}$ and any $q\in\mathbb{R}$ we have%
\[
\int_{p\geq q}\left|  \left\langle f|\mathbf{e}_{\lambda}\right\rangle
\right|  ^{2}\mathrm{d}\lambda=\int_{-\infty}^{\infty}I\left(  x-q\right)
\left|  f\left(  x\right)  \right|  ^{2}\mathrm{d}x
\]
where%
\[
I\left(  x\right)  \doteq2^{1/2}\int_{-\infty}^{x}\exp\left(  -2\pi
y^{2}\right)  \mathrm{d}y
\]
\end{corollary}

\textbf{Remark.} The function $I\left(  x\right)  $ tends fast to $1$ and $0$
as $x\rightarrow\infty,$ respectively, $x\rightarrow-\infty$ and $I\left(
0\right)  =1/2.$

\textsc{Proof }follows from (\ref{dt}). $\blacktriangleright$

\begin{corollary}
Let $D$ be a domain and $S$ be a rotation in the phase plane $\Phi$ and
$M_{S}$ be the corresponding metaplectic transform. We have
\begin{equation}
\int\left|  M_{S}f\left(  x\right)  \right|  ^{2}I\left(  x-q\right)
\mathrm{d}x\leq\int_{\Phi\backslash D}\left|  \left\langle f|\mathbf{e}%
_{\lambda}\right\rangle \right|  ^{2}\mathrm{d}\lambda\label{msf}%
\end{equation}
where $q=\sup\left\{  p;\lambda=\left(  p,\theta\right)  \in S\left(
D\right)  \right\}  .$
\end{corollary}

\textsc{Proof.} By Proposition \ref{hp} we have%
\[
\int\left|  M_{S}f\left(  x\right)  \right|  ^{2}I\left(  x-q\right)
\mathrm{d}x\leq\int_{\Phi\backslash S\left(  D\right)  }\left|  \left\langle
M_{S}f|\mathbf{e}_{\lambda}\right\rangle \right|  ^{2}\mathrm{d}\lambda.
\]
Because of the transform $M_{S}$ is unitary, we have by Corollary \ref{rot} we
have%
\[
\left\langle f|\mathbf{e}_{\lambda}\right\rangle =\left\langle M_{S}%
f|M_{S}\mathbf{e}_{\lambda}\right\rangle =\exp\left(  -\imath\phi\right)
\left\langle M_{S}f|\mathbf{e}_{S\left(  \lambda\right)  }\right\rangle
\]
for some real phase $\phi.$ Therefore
\[
\int_{\Phi\backslash S\left(  D\right)  }\left|  \left\langle M_{S}%
f|\mathbf{e}_{\lambda}\right\rangle \right|  ^{2}\mathrm{d}\lambda=\int
_{\Phi\backslash D}\left|  \left\langle M_{S}f|\mathbf{e}_{S\left(
\lambda\right)  }\right\rangle \right|  ^{2}\mathrm{d}\lambda=\int
_{\Phi\backslash D}\left|  \left\langle f|\mathbf{e}_{\lambda}\right\rangle
\right|  ^{2}\mathrm{d}\lambda
\]
which completes the proof. $\blacktriangleright$

Suppose that the right-hand side of (\ref{msf}) is small. Taking in account
that $I\left(  x\right)  \geq1/2$ for $x\geq0,$ we conclude that for any
rotation $S$ the function $M_{S}f$ is strongly concentrated in the interval
$\left[  r,q\right]  \doteq\pi\left(  S\left(  D\right)  \right)  ,$ where
$\pi:\Phi\rightarrow\mathbb{R}$ is the orthogonal projection. The inverse is
true if, for instance, $D$ is a convex polygon.

\section{Gabor series}

Choose some elements $\varepsilon\in\mathbb{R}$, $\varepsilon^{\ast}%
\in\mathbb{R}^{\prime}$ such that $\varepsilon^{\ast}\left(  \varepsilon
\right)  =1$ and consider the lattice $\Lambda\subset\Phi=\mathbb{R\times
R}^{\prime}$ generated by $\varepsilon$ and $\varepsilon^{\ast}$:%
\[
\Lambda\doteq\left\{  \lambda=k\varepsilon+j\varepsilon^{\ast},\;k,j\in
\mathbb{Z}\right\}  .
\]
The simplectic area of a cell is equal to $\varepsilon^{\ast}\left(
\varepsilon\right)  =1.$ Introduce a coordinate $x$ in $\mathbb{R}$ such that
$x\left(  \varepsilon\right)  =1;$ then we have $\xi\left(  \varepsilon^{\ast
}\right)  =1$ for he dual coordinate $\xi$ in $\mathbb{R}^{\prime}.$ Consider
the Gabor system $\left\{  \mathbf{e}_{\lambda},\,\lambda\in\Lambda\right\}  $.

\begin{proposition}
\label{p1} For an arbitrary numerical sequence $\left\{  c_{\lambda}\right\}
\in l_{2}\left(  \Lambda\right)  $ the series
\[
f=\sum_{\Lambda}c_{\lambda}\mathbf{e}_{\lambda}%
\]
converges in $L_{2}$ and the inequality holds
\[
\left\|  f\right\|  \leq\sigma_{0}\left(  \sum\left|  c_{\lambda}\right|
^{2}\right)  ^{1/2},\,\,\sigma_{0}=\sum\exp\left(  -\frac{\pi}{2}k^{2}\right)
.
\]
\end{proposition}

\textsc{Proof.}\textbf{\ }By (\ref{ee}) we have for any $\lambda,\mu\in
\Lambda,$ $\lambda=\left(  p,\theta\right)  $
\begin{align}
\left\langle \mathbf{e}_{\lambda}|\mathbf{e}_{\mu}\right\rangle  &
=g_{\lambda-\mu},\;\label{eeg}\\
g_{\lambda}  &  \doteq\exp\left(  \pi\imath p\theta-\frac{\pi}{2}\left(
p^{2}+\theta^{2}\right)  \right)  .\nonumber
\end{align}
Therefore
\[
\left\|  f\right\|  ^{2}=\sum_{\lambda,\mu}c_{\lambda}\bar{c}_{\mu
}\left\langle \mathbf{e}_{\lambda}|\mathbf{e}_{\mu}\right\rangle
=\sum_{\lambda}c_{\lambda}\sum_{\mu}\bar{c}_{\mu}g_{\lambda-\mu}%
\]
The interior sum is a convolution and we can estimate its $l_{2}$-norm as
follows:
\[
\sum_{\lambda}\left|  \sum_{\mu}\bar{c}_{\mu}g_{\lambda-\mu}\right|  ^{2}%
\leq\left(  \sum\left|  g_{\lambda}\right|  \right)  ^{2}\sum\left|  c_{\mu
}\right|  ^{2}=\sigma_{0}^{4}\sum\left|  c_{\lambda}\right|  ^{2},
\]
since $\sum\left|  g_{\lambda}\right|  =\sigma_{0}^{2}.$ This yields
\begin{equation}
\left\|  f\right\|  ^{4}\leq\sum_{\lambda}\left|  c_{\lambda}\right|  ^{2}%
\sum_{\lambda}\left|  \sum_{\mu}\bar{c}_{\mu}g_{\lambda-\mu}\right|  ^{2}%
\leq\sigma_{0}^{4}\left(  \sum\left|  c_{\mu}\right|  ^{2}\right)
^{2}.\,\blacktriangleright\label{escon}%
\end{equation}

\section{Zak transform}

The Zak transform of a function $f\in L_{2}$ is defined by the series
\[
Zf\left(  y,\xi\right)  =\sum_{q\in\mathbb{Z}}\exp\left(  2\pi\imath
q\xi\right)  f\left(  y+q\right)  ,
\]
which converges almost everywhere in the square $\mathbb{Q}\doteq\left\{
0\leq y,\xi\leq1\right\}  .$ The function $Zf$ fulfils
\begin{equation}
Zf\left(  y,\xi+1\right)  =Zf\left(  y,\xi\right)  ,\;Zf\left(  y+1,\xi
\right)  =\exp\left(  -2\pi\imath\xi\right)  Zf\left(  y,\xi\right)  .
\label{zz}%
\end{equation}

\begin{proposition}
The Zak transform is a unitary operator $L_{2}\left(  \mathbb{R}\right)
\rightarrow L_{2}\left(  \mathbb{Q}\right)  .$
\end{proposition}

\textsc{Proof} is straightforward and we omit it. The inversion formula reads%
\begin{align*}
f\left(  x\right)   &  =\sum_{r\in\mathbb{Z}}\int_{0}^{1}\int_{0}%
^{1}\mathrm{d}y\mathrm{d}\xi\,g\left(  y,\xi\right)  \exp\left(  2\pi
\imath\left(  y-x\right)  \left(  \xi+r\right)  \right) \\
&  =\int_{0}^{1}\mathrm{d}y\int_{\mathbb{R}}\,g\left(  y,\xi\right)
\exp\left(  2\pi\imath\xi\left(  y-x\right)  \right)  \mathrm{d}\xi,
\end{align*}
where $g=Zf$ is 1-periodic function of $\xi.$

The Weyl-Heisenberg group acts in $L_{2}$ by the shifts and modulation
operators
\begin{equation}
T_{\lambda}f\left(  x\right)  =\exp\left(  2\pi\imath\theta x\right)  f\left(
x-p\right)  ,\lambda=\left(  p,\theta\right)  \label{WH}%
\end{equation}
We have for any $\lambda\in\Lambda$ we have%
\begin{equation}
Z\left(  T_{\lambda}f\right)  \left(  y,\xi\right)  =\exp\left(  2\pi
\imath\left(  p\xi+\theta\xi\right)  \right)  Zf\left(  y,\xi\right)  .
\label{ztm}%
\end{equation}
The sum
\[
\Theta\left(  z\right)  =2^{1/4}\sum_{q\in\mathbb{Z}}\exp\left(  2\pi\imath
qz-\pi q^{2}\right)
\]
is a Jacobi elliptic function: $\Theta\left(  z\right)  =2^{1/2}\theta
_{3}\left(  z;\imath\right)  $. This function\ is holomorphic of $z$ in the
whole plane, satisfies the periodicity conditions
\begin{equation}
\Theta\left(  z+1\right)  =\Theta\left(  z\right)  ,\;\Theta\left(
z+\imath\right)  =\exp\left(  \pi-2\pi\imath z\right)  \Theta\left(  z\right)
\label{tt}%
\end{equation}
has simple zero at the point $1/2+\imath/2$ and no other zeros in the closed
square $\mathbb{Q}$. The Zak transform of Gabor function for $\lambda=\left(
p,\theta\right)  \in\Lambda$ is expressed in terms of $\Theta:$
\begin{equation}
Z\mathbf{e}_{\lambda}\left(  y,\xi\right)  =\exp\left(  2\pi\imath\left(
p\xi+\theta y\right)  \right)  \,\exp\left(  -\pi y^{2}\right)  \Theta\left(
\xi+\imath y\right)  . \label{ze}%
\end{equation}

\section{Creation and annihilation operators}

The operator
\[
\mathbf{a}=\frac{1}{2\pi}\frac{\mathrm{d}}{\mathrm{d}x}+x,\; \mathbf{a}%
^{+}=-\frac{1}{2\pi}\frac{\mathrm{d}}{\mathrm{d}x}+x
\]
in $L_{2}$ are adjoint one to another. They are called the
\textit{annihilation} and the \textit{creation} operators. Any Gabor function
is an eigenvector of the annihilation operator:
\begin{equation}
\mathbf{ae}_{\lambda}=\mathbf{\lambda e}_{\lambda} \label{ei}%
\end{equation}
where $\lambda=\left(  p,\theta\right)  $ and $\,\mathbf{\lambda}%
=p+\imath\theta$. For any $\phi$ in the domain of the operator $\mathbf{a}$ we
have
\begin{equation}
Z\left(  \mathbf{a}\phi\right)  =AZ\phi,\,A\doteq\frac{1}{2\pi\imath}\left(
\frac{\partial}{\partial\xi}+\imath\frac{\partial}{\partial y}\right)  +y.
\label{ZZ}%
\end{equation}

Consider the line bundle $\mathbf{E}$ over the torus $\mathbb{R}^{2}/\Lambda$
whose sections $s$ fulfil the periodicity conditions%
\begin{equation}
s\left(  y,\xi+1\right)  =s\left(  y,\xi\right)  ,\;s\left(  y+1,\xi\right)
=\exp\left(  -2\pi\imath\xi\right)  s\left(  y,\xi\right)  . \label{ss}%
\end{equation}
in the unit square $\mathbb{Q}$ (which is the fundamental domain of the
lattice $\Lambda$). Let $W_{2}^{\alpha}\left(  \mathbf{E}\right)  $ be the
Sobolev space of order $\alpha\geq0$ of sections $\phi$ of $\mathbf{E.}$ By
(\ref{zz}) the function $Zf$ is a section of $\mathbf{E}$ of the class
$W_{2}^{0}\left(  \mathbf{E}\right)  =L_{2}\left(  \mathbb{Q}\right)  .$ The
differential operator $A$ acts on smooth sections of $\mathbf{E}$ and defines
for any $\alpha$ a bounded map $A:W^{\alpha+1}\left(  \mathbf{E}\right)
\rightarrow W^{\alpha}\left(  \mathbf{E}\right)  .$

\begin{proposition}
\label{w}For any function $f\in L_{2}$ such that%
\begin{equation}
\int\left(  \left|  \lambda\right|  ^{\delta}+1\right)  \left|  \left\langle
f|\mathbf{e}_{\lambda}\right\rangle \right|  ^{2}\mathrm{d}\lambda
<\infty,\;\delta>0 \label{lfe}%
\end{equation}
we have $Zf\in W_{2}^{\delta}(\mathbf{E}).$
\end{proposition}

\textsc{Proof.} By (\ref{fed}) we can write%
\[
Zf=\int_{\Phi}\left\langle f|\mathbf{e}_{\lambda}\right\rangle Z\mathbf{e}%
_{\lambda}\mathrm{d}\lambda
\]
We have for arbitrary $\lambda=\left(  p,\theta\right)  \in\Phi$
\[
Z\mathbf{e}_{\lambda}\left(  \xi,y\right)  =2^{1/4}\sum_{q\in\mathbb{Z}}%
\exp\left(  2\pi\imath q\xi+\theta\left(  y+q\right)  \right)  \exp\left(
-\pi\left(  y+q-p\right)  ^{2}\right)
\]
Setting $p=k+r,\theta=j+\eta,\,k,l\in\mathbb{Z}$ and $q^{\prime}=q-k$ yields%
\begin{align}
Z\mathbf{e}_{\lambda}\left(  y,\xi\right)   &  =2^{1/4}\sum_{q\in\mathbb{Z}%
}\exp\left(  2\pi\imath\left[  q\xi+\theta y+q\eta\right]  \right)
\exp\left(  -\pi\left(  y+q-k-r\right)  ^{2}\right) \nonumber\\
&  =2^{1/4}\exp\left(  2\pi\imath\left[  k\left(  \xi+\eta\right)  +\theta
y\right]  \right)  \sum_{q\in\mathbb{Z}}\exp\left(  2\pi\imath q^{\prime
}\left(  \xi+\eta\right)  \right)  \exp\left(  -\pi\left(  y+q^{\prime
}-r\right)  ^{2}\right) \nonumber\\
&  =\exp\left(  2\pi\imath\left[  k\xi+jy\right]  \right)  \exp\left(
2\pi\imath k\eta\right)  Z\mathbf{e}_{\left(  r,\eta\right)  }\left(
y,\xi\right)  \label{zel}%
\end{align}
We can write%
\[
Zf\left(  y,\xi\right)  =\sum_{\mathbb{Z\times Z}}\exp\left(  2\pi
\imath\left[  k\xi+jy\right]  \right)  \int_{\mathbb{Q}}\left\langle
f|\mathbf{e}_{\lambda}\right\rangle E_{\lambda}\left(  y,\xi\right)
\mathrm{d}r\,\mathrm{d}\eta,
\]
The kernel $E_{\lambda}\left(  y,\xi\right)  \doteq\exp\left(  2\pi\imath
k\eta\right)  Z\mathbf{e}_{\left(  r,\eta\right)  }\left(  y,\xi\right)  $ is
a real analytic function for $\lambda,\left(  y,\xi\right)  \in\mathbb{R\times
R}$ , is 1-periodic in $\xi$ and its derivatives in $y,\xi$ are bounded in any
strip $\left|  y\right|  \leq C$ uniformly for $\lambda\in\Phi.$ Take a test
function $\phi$ in $\mathbb{R\times R}$ and apply the Fourier transform
$F=F_{\left(  y,\xi\right)  \mapsto\left(  s,t\right)  }$ :%
\[
F\left(  \phi Zf\right)  \left(  s,t\right)  =\int_{\mathbb{Q}}\sum
_{k,j}\left\langle f|\mathbf{e}_{\lambda}\right\rangle F\left(  \phi
E_{\lambda}\right)  \left(  s-k,t-j\right)  \mathrm{d}r\,\mathrm{d}\eta
\]
The interior sum can be written as convolution on $\mathbb{Z}^{2}$:%
\[
G_{r,\eta}\doteq\sum_{k,j\in\mathbb{Z}}\left\langle f|\mathbf{e}_{\left(
k,j\right)  +\left(  r,\eta\right)  }\right\rangle F\left(  \phi E_{\lambda
}\right)  \left(  s-k,t-j\right)
\]
We have%
\[
\left|  F\left(  \phi E_{\lambda}\right)  \left(  s,t\right)  \right|  \leq
C_{n}\left(  \left|  s\right|  +\left|  t\right|  +1\right)  ^{-n}%
\]
for any natural $n,$ since $\phi E_{\lambda}\in C^{\infty}.$ Therefore for any
fixed $\left(  r,\eta\right)  \in\mathbb{Q}$%
\[
\int\left(  \left|  s\right|  +\left|  t\right|  +1\right)  ^{\delta}\left|
G_{r,\eta}\left(  s,t\right)  \right|  ^{2}\mathrm{d}s\,\mathrm{d}t\leq
C\sum_{k,j}\left(  \left|  k\right|  +\left|  j\right|  +1\right)  ^{\delta
}\left|  \left\langle f|\mathbf{e}_{\left(  k,j\right)  +\left(
r,\eta\right)  }\right\rangle \right|  ^{2}.
\]
where the sum in the right-hand side converges for almost all $\left(
r,\eta\right)  $ and the constant $C$ does not depend on $f.$ Integrating over
$\mathbb{Q}$ yields%
\[
\int\left(  \left|  s\right|  +\left|  t\right|  +1\right)  ^{\delta}\left|
F\left(  \phi Zf\right)  \left(  s,t\right)  \right|  ^{2}\mathrm{d}%
s\,\mathrm{d}t= \int\left(  \left|  s\right|  +\left|  t\right|  +1\right)
^{\delta}\left|  \int_{\mathbb{Q}}G_{r,\eta}\left(  s,t\right)  \mathrm{d}%
r\,\mathrm{d}\eta\right|  ^{2}\mathrm{d}s\mathrm{d}t
\]
\begin{align*}
&  \leq\int_{\mathbb{Q}}\int\left(  \left|  s\right|  +\left|  t\right|
+1\right)  ^{\delta}\left|  G_{r,\eta}\left(  s,t\right)  \right|
^{2}\mathrm{d}s\,\mathrm{d}t\,\mathrm{d}r\,\mathrm{d}\eta\\
&  \leq C\int_{\mathbb{Q}}\sum_{k,j}\left(  \left|  k\right|  +\left|
j\right|  +1\right)  ^{\delta}\left|  \left\langle f|\mathbf{e}_{\left(
k,j\right)  +\left(  r,\eta\right)  }\right\rangle \right|  ^{2}%
\mathrm{d}r\,\mathrm{d}\eta\\
&  =C\int_{\Phi}\left(  \left|  k\right|  +\left|  j\right|  +1\right)
^{\delta}\left|  \left\langle f|\mathbf{e}_{\lambda}\right\rangle \right|
^{2}\mathrm{d}p\,\mathrm{d}\theta\\
&  \leq C^{\prime}\int_{\Phi}\left(  \left|  \lambda\right|  +1\right)
^{\delta}\left|  \left\langle f|\mathbf{e}_{\lambda}\right\rangle \right|
^{2}\mathrm{d}\lambda<\infty
\end{align*}
This implies that $\phi Zf$ belongs to $W_{2}^{\delta}$ and our statement
follows. $\blacktriangleright$

\textbf{Definition.} Denote by $H^{\delta}$ the space of functions $f\in
L_{2}$ that satisfies the condition (\ref{lfe}) and set
\[
\left\|  f\right\| _{\delta}=\left(  \int\left(  \left|  \lambda\right|
^{\delta}+1\right)  \left|  \left\langle f|\mathbf{e}_{\lambda}\right\rangle
\right|  ^{2}\mathrm{d}\lambda\right)  ^{1/2}%
\]
Denote$\;H\doteq\cup_{\delta>1}H^{\delta}.$

\textbf{Remark 1.} The space $H^{\delta}$ coincides with the modulation spaces
$\mathbf{M}_{2,2}^{w}$ of H.Feichtinger, \cite{FG} for the weight function
$w\left(  \lambda\right)  =\left|  \lambda\right|  ^{\delta}+1.$

\textbf{Remark 2.} The domain of the harmonic oscillator operator
\[
\mathrm{H}\doteq\frac{1}{2}\left(  \mathbf{a}^{+}\mathbf{a}+\mathbf{aa}%
^{+}\right)  =-\frac{1}{4\pi^{2}}\frac{\mathrm{d}^{2}}{\mathrm{d}x^{2}}%
+x^{2}.
\]
is equal to the space $H^{2}.$

\section{Relaxed Gabor expansion}

We show now that an arbitrary function $f\in H$ can be expanded in a Gabor
series as above with one more term.

Denote $\sharp=\left(  1/2,1/2\right)  \in\Phi;$ we call \textit{sharp} point
any element $\mu\in\Lambda+$ $\sharp$, \footnote{The motivation of this term
comes from the interpretation of the lattice $\Lambda$ as a table for a
regular time-frequency notation of music (Brillouin-Wigner). The point
$\sharp$ is then a tune which is shifted inside a rhythm interval (syncope),
whose height is shifted by half a tune (sharp= ''diese'' tune).}. The
\textit{sharp Poisson }functional is the series
\[
\gamma^{\sharp}\left(  f\right)  \doteq\frac{1}{\imath\Theta\left(  0\right)
}\sum_{\mathbb{Z}}\left(  -1\right)  ^{q}f\left(  q+1/2\right)  =\frac
{Zf\left(  \sharp\right)  }{\imath\Theta\left(  0\right)  }.
\]
We show below that this functional is well defined and continuous in $H$. In
particular, we have $\gamma^{\sharp}\left(  \mathbf{e}_{\sharp}\right)
=1,$whereas $\gamma^{\sharp}\left(  \mathbf{e}_{\lambda}\right)  =0$ for
arbitrary $\lambda\in\Lambda,$ since $\Theta$ vanishes at each sharp point. We
call the set $\Lambda^{\sharp}\doteq\Lambda\cup\left\{  \sharp\right\}  $ the
\textit{relaxed} lattice and consider the relaxed\textit{ }Gabor system
$\left\{  \mathbf{e}_{\lambda},\,\lambda\in\Lambda^{\sharp}\right\}  $.

\begin{theorem}
\label{t1}There exists a family of continuous functionals $\gamma^{\lambda
},\lambda\in\Lambda$ in the space $H$ such that an arbitrary $f\in H$ is
developed in the series
\begin{equation}
f=\sum_{\Lambda^{\sharp}}\gamma^{\lambda}\left(  f\right)  \mathbf{e}%
_{\lambda} \label{ge}%
\end{equation}
that converges in $L_{2}\left(  \mathbb{R}\right)  $ and for any $\delta>1$
there exists a constant $C_{\delta}$ such that
\begin{equation}
\sum_{\Lambda^{\sharp}}\left|  \gamma^{\lambda}\left(  f\right)  \right|
^{2}\leq C_{\delta}\left\|  f\right\|  _{\delta}^{2}. \label{gd}%
\end{equation}
\end{theorem}

\begin{lemma}
\label{l1}For an arbitrary $f\in H$ the Zak transform $Zf$ is $\varepsilon
$-H\"{o}lder continuous for $\varepsilon<\delta-1.$
\end{lemma}

\textsc{Proof of Lemma.} By Proposition \ref{w} $Zf\in W_{2}^{\delta}\left(
\mathbf{E}\right)  $ for some $\delta>1.$ Sobolev's imbedding theorem implies
that $Zf$ is a continuous section of $\mathbf{E}$ and $\sup\left|  Zf\right|
\leq C\left\|  f\right\|  _{\delta},$ \cite{Sm}. Moreover, $Zf$ belongs to the
H\"{o}lder $\varepsilon$-class for $\varepsilon<\delta-1$.
$\blacktriangleright$

\textsc{Proof of Theorem.} By Lemma \ref{l1} the functional $\gamma^{\sharp}$
is well defined in $H.$ Set $f_{\sharp}=f-\gamma^{\sharp}\left(  f\right)
\mathbf{e}_{\sharp}\in H$ and define $F\doteq Zf_{\sharp}/Z\mathbf{e}_{0}.$ By
Lemma \ref{l3} the function $F$ is square integrable. It is double periodic
with unit periods and can be represented by the double Fourier series
\[
F\left(  y,\xi\right)  =\sum_{\lambda=\left(  p,\theta\right)  \in\Lambda
}c_{\lambda}\exp\left(  2\pi\imath\left(  py+\theta\xi\right)  \right)  ,
\]
where $\sum\left|  c_{\lambda}\right|  ^{2}=\left\|  F\right\|  ^{2}.$ By
Proposition \ref{p1} the series $\sum c_{\lambda}\mathbf{e}_{\lambda}$
converges in $L_{2}\left(  \mathbb{R}\right)  $ to a function $g$. On the
other hand, because of (\ref{ze}),
\[
Zf_{\sharp}\left(  y,\xi\right)  =Z\mathbf{e}_{0}\sum c_{\lambda}\exp\left(
2\pi\imath\left(  py+\theta\xi\right)  \right)  =\sum c_{\lambda}%
Z\mathbf{e}_{\lambda}\left(  y,\xi\right)
\]
This implies $f_{\sharp}=g,$ since $Z$ is a unitary operator. Set
$\gamma^{\lambda}\left(  f\right)  \doteq c_{\lambda}$ for $\lambda\in\Lambda
$. $\blacktriangleright$

\begin{lemma}
\label{l3}We have $F\in L_{2}\left(  \mathbb{Q}\right)  .$
\end{lemma}

\textsc{Proof of Lemma}. By the previous Lemma the function $Zf_{\sharp}$ is
$\varepsilon$-H\"{o}lder continuous in the interior of $\mathbb{Q}.$ It
vanishes at the sharp point, since $\gamma^{\sharp}\left(  f_{\sharp}\right)
=0.$ Therefore $\left|  Zf_{\sharp}\left(  z\right)  \right|  \leq C_{\delta
}\left|  z-1/2-\imath/2\right|  ^{\varepsilon}.$ This implies
\[
\left|  F\left(  z\right)  \right|  =\left|  \frac{\exp\left(  \pi
y^{2}\right)  Zf_{\sharp}}{\Theta\left(  z\right)  }\right|  \leq C\left|
z-1/2-\imath/2\right|  ^{\varepsilon-1},
\]
since $\Theta\left(  z\right)  $ has simple zero at each sharp point. This
implies the statement. $\blacktriangleright$

\textbf{Remark 1. }More strong inequality holds for $f\in H^{\delta}%
,\delta>1:$%
\[
\sum_{\Lambda}\left(  \left|  \lambda\right|  +1\right)  ^{2\varepsilon
}\left|  \gamma^{\lambda}\left(  f\right)  \right|  ^{2}\leq C_{\delta
,\varepsilon}\left\|  f\right\|  _{\delta}^{2},\;0<\varepsilon<\delta-1
\]

\textbf{Remark 2}. We can take an arbitrary sharp point $\mu\in\Lambda+\sharp$
instead of $\left(  1/2.1/2\right)  $ in Theorem \ref{t1}.

\section{Gabor coefficients}

The coefficients are uniquely defined since of

\begin{proposition}
\label{u} If $\left\{  c_{\lambda}\right\}  \in l_{2}$ $\left(  \Lambda
\right)  ,b\in\mathbb{C}$ and
\begin{equation}
b\mathbf{e}_{\sharp}+\sum_{\Lambda}c_{\lambda}\mathbf{e}_{\lambda}=0,
\label{sz}%
\end{equation}
then $c_{\lambda}=0$ for all $\lambda$.
\end{proposition}

\textsc{Proof.} By Proposition \ref{p1} the series (\ref{sz}) converges in
$L_{2}$. Apply the Zak transform. By (\ref{ze}) we get
\[
-bZ\mathbf{e}_{\sharp}=\sum c_{\lambda}Z\mathbf{e}_{\lambda}=Z\mathbf{e}%
_{0}\sum c_{\lambda}\exp2\pi\imath\left(  p\xi+\theta y\right)
\]
The series in the right-hand side converges to a function $g\in L_{2}\left(
\mathbb{Q}\right)  $ since $Z$ is a unitary operator. It follows
$g=-bZ\mathbf{e}_{\sharp}/Z\mathbf{e}_{0}.$ We have $Z\mathbf{e}_{\sharp
}\left(  1/2,1/2\right)  \neq0,$ whereas the function $Z\mathbf{e}_{0}\left(
y,\xi\right)  =\exp\left(  -\pi y^{2}\right)  \Theta\left(  \xi+\imath
y\right)  $ vanishes in the point $1/2+\imath/2.$ Therefore the function $g$
can not be square integrable unless $b=0.$ This yields $g=0,$ hence
$c_{\lambda}=0$ for all $\lambda\in\Lambda.$ $\blacktriangleright$

\begin{corollary}
We have the following explicit formula for the coefficients:
\begin{equation}
\gamma_{\lambda}\left(  f\right)  =\int_{0}^{1}\int_{0}^{1}\exp\left(
-2\pi\imath\left(  py+\theta\xi\right)  \right)  \frac{\exp\left(  \pi
y^{2}\right)  Zf_{\sharp}\left(  y,\xi\right)  }{\Theta\left(  \xi+\imath
y\right)  }\mathrm{d}y\mathrm{d}\xi,\;\lambda\in\Lambda\label{cla}%
\end{equation}
where $f_{\sharp}=f-\gamma^{\sharp}\left(  f\right)  \mathbf{e}_{\sharp}.$
\end{corollary}

\textbf{Remark. }It follows from a formula of M.Bastiaans \cite{Ba}, that the
coefficients can be written in the form
\[
\gamma^{\lambda}\left(  f\right)  =\left\langle f_{\sharp}|\exp\left(
2\pi\imath\theta x\right)  \gamma\left(  x-p\right)  \right\rangle
,\;\lambda=\left(  p,\theta\right)  \in\Lambda,
\]
where
\[
\gamma\left(  x\right)  =2^{-1/4}\left(  \pi/K_{0}\right)  ^{3/2}\exp\left(
\pi x^{2}\right)  \sum_{n+1/2\geq x}\left(  -1\right)  ^{n}\exp\left(
-\pi\left(  n+1/2\right)  ^{2}\right)
\]
This function $\gamma$ does not belong to $L_{2},$ but the integrals
$\left\langle f_{\sharp}|\exp\left(  2\pi\imath\theta x\right)  \gamma\left(
x-p\right)  \right\rangle $ converge since $f_{\sharp}\in H$ and
$\gamma^{\sharp}\left(  f_{\sharp}\right)  =0.$

\section{Improving convergence}

To ensure faster convergence of the relaxed Gabor series we impose more sharp
conditions on $f.$ Fix an integer $m>0$ and consider the norm
\[
\left\|  f\right\|  _{\delta,m}\doteq\left(  \sum_{j=0}^{m}\left\|
\mathbf{a}^{j}f\right\|  _{\delta}^{2}\right)  ^{1/2}%
\]
Let $H^{\delta,m}$ be the space of functions with finite norm $\left\|
\cdot\right\|  _{\delta,m}.$

\begin{lemma}
\label{van}For arbitrary different points $\mu_{0},...,\mu_{m}\in\mathbb{C}$
the inverse to the VanderMond matrix $W=\{w_{j}^{k}\doteq(\mu_{j})^{k}\}$ is
equal to the matrix $V=\{v_{j}^{k}\}$ where
\begin{align*}
v_{j}^{k}  &  \doteq\frac{\sigma_{m+1-j}^{\left(  k\right)  }}{p^{\prime
}\left(  \mu_{k}\right)  },\,p\left(  \lambda\right)  =\Pi_{j=0}^{m}\left(
\lambda-\mu_{j}\right)  ,\\
\sigma_{j}^{\left(  k\right)  }  &  =\left(  -1\right)  ^{j}\sigma_{j}\left(
\mu_{0},...,\widehat{\mu_{k}},...,\mu_{m}\right)
\end{align*}
and $\sigma_{j}$ denotes the $j$-th elementary symmetric polynomial.
\end{lemma}

\textsc{Proof} is by straightforward check.

\begin{lemma}
\label{det}For an arbitrary natural $m$ and any set of different points
$\mu_{0},...,\mu_{m}\in\Lambda+\sharp$ the are numbers$\left\{  h_{j}%
^{k}\right\}  $ such that%
\begin{equation}
\gamma^{\sharp}\left(  \mathbf{a}^{k}\mathbf{d}_{m,j}\right)  =\delta_{j}%
^{k},\,k,j=0,...,m \label{gdd}%
\end{equation}
where%
\begin{equation}
\mathbf{d}_{m,j}=\sum_{s=0}^{m}h_{j}^{s}\mathbf{e}_{\mu_{s}},\;j=0,...,m.
\label{dhe}%
\end{equation}
\end{lemma}

\textsc{Proof.} Calculate the matrix
\[
g_{j}^{k}\doteq\gamma^{\sharp}\left(  \mathbf{a}^{k}\mathbf{e}_{\mu_{j}%
}\right)  ,\;j,k=0,...,m
\]
By (\ref{zel}) we have for an arbitrary $\mu=\left(  q,\eta\right)  \in
\Lambda+\sharp$%
\begin{align*}
Z\mathbf{e}_{\mu}\left(  y,\xi\right)   &  =\exp\left(  -\pi\left(
y-1/2\right)  ^{2}+2\pi\imath\left[  y\eta+\left(  q-1/2\right)  \left(
\xi+1/2\right)  \right]  \right) \\
&  \times\Theta\left(  \xi+1/2+\imath\left(  y-1/2\right)  \right)
\end{align*}
which yields $Z\mathbf{e}_{\mu}\left(  \sharp\right)  =\exp\left(  \pi
\imath\eta\right)  \Theta\left(  0\right)  .$ By (\ref{ei})
\[
\gamma^{\sharp}\left(  \mathbf{a}^{k}\mathbf{e}_{\mu}\right)  =\frac
{\mathbf{\mu}^{k}}{\imath\Theta\left(  0\right)  }Z\mathbf{e}_{\mu}\left(
\sharp\right)  =\left(  -1\right)  ^{\left[  \eta\right]  }\mathbf{\mu}^{k}%
\]
This yields $g_{j}^{k}=\left(  -1\right)  ^{\left[  \eta_{j}\right]
}\mathbf{\mu}_{j}^{k}.$ The matrix $\left\{  g_{j}^{k}\right\}  $ is
invertible and by Lemma \ref{van} the entries of the inverse matrix are%
\[
h_{j}^{k}=\left(  -1\right)  ^{\left[  \eta_{k}\right]  }\frac{\sigma
_{m+1-j}^{\left(  k\right)  }}{p^{\prime}\left(  \mathbf{\mu}_{k}\right)
},\,k,j=0,...,m.
\]
Then (\ref{dhe}) implies (\ref{gdd}). $\blacktriangleright$

\begin{lemma}
\label{d}We have for any $\lambda\in\Phi$%
\begin{equation}
\left\|  \sum\mathbf{d}_{m,j}\mathbf{\lambda}^{j}\right\|  _{\delta,m}%
\leq\left(  M+1\right)  ^{m+\delta}L^{m} \label{dml}%
\end{equation}
where
\[
M\doteq\max\left|  \mu_{s}\right|  ,\;L=\max\left|  \lambda-\mu_{s}\right|
\]
\end{lemma}

\textsc{Proof.} We have for $k,j=1,...,m$%
\[
\left\|  \mathbf{a}^{k}\sum_{j}\mathbf{d}_{m,j}\mathbf{\lambda}^{j}\right\|
=\left\|  \sum_{j,s}h_{j}^{s}\mathbf{\lambda}^{j}\mathbf{\mu}_{s}%
^{k}\mathbf{e}_{\mu_{s}}\right\|  \leq M^{k}\left|  \sum_{j}h_{j}%
^{s}\mathbf{\lambda}^{j}\right|
\]
since $\left\|  \mathbf{e}_{\mu}\right\|  =1.$ Further%
\[
\sum_{j}h_{j}^{s}\mathbf{\lambda}^{j}=\pm\sum_{j}\frac{\sigma_{m+1-j}^{\left(
s\right)  }}{p^{\prime}\left(  \mathbf{\mu}_{s}\right)  }\mathbf{\lambda}%
^{j}=\pm\frac{p_{s}\left(  \mathbf{\lambda}\right)  }{p_{s}\left(
\mathbf{\mu}_{s}\right)  }%
\]
where $p_{s}\left(  \mathbf{\lambda}\right)  \doteq\Pi_{j\neq s}\left(
\mathbf{\lambda}-\mathbf{\mu}_{j}\right)  $. This yields%
\[
\left|  \sum_{j}h_{j}^{s}\mathbf{\lambda}^{j}\right|  =\left|  \frac
{p_{s}\left(  \mathbf{\lambda}\right)  }{p_{s}\left(  \mathbf{\mu}_{s}\right)
}\right|  =\prod_{j\neq s}\left|  \frac{\mathbf{\lambda-\mu}_{j}}{\mathbf{\mu
}_{s}\mathbf{-\mu}_{j}}\right|  \leq L^{m}%
\]
and (\ref{dml}) follows. $\blacktriangleright$

\begin{theorem}
\label{m}For an arbitrary natural $m$ and different points $\mu_{0}%
,...,\mu_{m}\in\Lambda+\sharp$ there exists a family of continuous functionals
$\gamma_{m}^{\lambda},$ $\lambda\in\Lambda$\ in $H^{m}=\cup_{\delta
>1}H^{\delta,m}$ such that for an arbitrary $f\in H^{m}$ the equation holds in
the space $W_{2}^{k}\left(  \mathbb{R}\right)  :$
\begin{equation}
f=\sum_{j=0}^{m}\gamma^{\sharp}\left(  \mathbf{a}^{j}f\right)  \mathbf{d}%
_{m,j}+\sum_{\Lambda}\gamma_{m}^{\lambda}\left(  f\right)  \mathbf{e}%
_{\lambda}, \label{keq}%
\end{equation}
where $\mathbf{d}_{m,j}$ are as in (\ref{dhe}) and there exists for any
$\delta>1$ a constant $C_{\delta}$ such that
\begin{equation}
\sum_{\Lambda}\left(  \left|  \lambda\right|  ^{2}+1\right)  ^{m}\left|
\gamma_{m}^{\lambda}\left(  f\right)  \right|  ^{2}\leq C_{\delta}\left\|
f\right\|  _{\delta,m}^{2}. \label{lgf}%
\end{equation}
\end{theorem}

\textsc{Proof.} Set
\[
f_{\sharp}=f-\sum_{j=0}^{m}\gamma^{\sharp}\left(  \mathbf{a}^{j}f\right)
\mathbf{d}_{j}\in H^{m}%
\]
By (\ref{ZZ}) and Lemma \ref{l1} the function $A^{j}Zf_{\sharp}$ is continuous
in the interior of $\mathbb{Q}$ for $j\leq m.$ We have $A=\bar{\partial}%
/2\pi\imath+y$, where $\bar{\partial}=\left(  \partial_{\xi}+\imath
\partial_{y}\right)  $ and
\[
A^{j}Zf_{\sharp}\left(  y,\xi\right)  =\Theta\left(  \xi+\imath y\right)
A^{j}\left[  \exp\left(  -\pi y^{2}\right)  F\left(  y,\xi\right)  \right]
=\Theta\left(  \xi+\imath y\right)  \exp\left(  -\pi y^{2}\right)
\bar{\partial}^{j}F\left(  y,\xi\right)  ,
\]
since the function $\Theta$ is holomorphic and $A\exp\left(  -\pi
y^{2}\right)  =0$. Therefore the function $\bar{\partial}^{j}F$ is H\"{o}lder
continuous in $\mathbb{Q}$ for $j\leq m.$ The property (\ref{gdd}) and
(\ref{ZZ}) imply that $\gamma^{\sharp}\left(  A^{j}Zf_{\sharp}\right)  =0,$
$j=0,...,m$ that is the left-hand side vanishes at the sharp point. By Lemma
\ref{l3} the function
\[
\bar{\partial}^{j}F=\frac{\exp\left(  \pi y^{2}\right)  A^{j}Zf_{\sharp
}\left(  y,\xi\right)  }{\Theta\left(  \xi+\imath y\right)  }%
\]
is double periodic and belongs to $L_{2}\left(  \mathbb{Q}\right)  $ for
$j=0,...,m.$ Therefore $F\in W_{2}^{m}\left(  \mathbb{R}^{2}/\Lambda\right)  $
and the Fourier coefficients $c_{\lambda}$ of $F$ satisfy
\[
\sum_{\Lambda}\left(  \left|  \lambda\right|  ^{2}+1\right)  ^{m}\left|
c_{\lambda}\right|  ^{2}\leq C\left\|  F\right\|  _{W_{2}^{m}}^{2}\leq
C^{\prime\prime}\left\|  f_{\sharp}\right\|  _{\delta,m}^{2}\leq
C^{\prime\prime\prime}\left\|  f\right\|  _{\delta,m}^{2}%
\]
This yields (\ref{lgf}) for $\gamma_{m}^{\lambda}\left(  f\right)  \doteq
c_{\lambda},$ $\lambda\in\Lambda.$ $\blacktriangleright$

A representation like (\ref{keq}) is unique, in spite of the additional terms:

\begin{proposition}
If for some integer $m\geq0$ the series
\begin{equation}
\sum_{\Lambda}\left(  \left|  \lambda\right|  ^{2}+1\right)  ^{m}\left|
c_{\lambda}\right|  ^{2} \label{lk}%
\end{equation}
converges and
\begin{equation}
\sum_{j=0}^{m}b_{j}\mathbf{d}_{j}+\sum_{\Lambda}c_{\lambda}\mathbf{e}%
_{\lambda}=0 \label{bce}%
\end{equation}
for some $b_{0},...,b_{m},$ then $b_{0}=...=b_{m}=0$ and $c_{\lambda}=0$ for
all $\lambda$.
\end{proposition}

\textsc{Proof.} The series in (\ref{bce}) converges to a function $g\in
W_{2}^{m}\left(  \mathbb{R}\right)  .$ Moreover $g$ belongs to the domain of
the operators $\mathbf{a}^{j},j=1,...,m.$ This follows from convergence of
(\ref{lk}) and arguments of Proposition \ref{p1}. By (\ref{ZZ}) the Zak
transform $Zf$ belongs to the domain of operators $A^{j},j=1,...,m.$ By
(\ref{ze}) we have%
\[
-\sum_{0}^{m}b_{j}Z\mathbf{d}_{j}=\sum_{\Lambda}c_{\lambda}Z\mathbf{e}%
_{\lambda}=Z\mathbf{e}_{0}\sum_{\Lambda}c_{\lambda}\exp2\pi\imath\left(
p\xi+\theta y\right)
\]
The series in the right-hand side converges to a function $g\in L_{2}\left(
\mathbb{Q}\right)  $ and $g=-\sum b_{j}Z\mathbf{d}_{j}/Z\mathbf{e}_{0}.$ We
have $Z\mathbf{d}_{0}\left(  z_{\sharp}\right)  =1,$ whereas the functions
$Z\mathbf{d}_{1},...,Z\mathbf{d}_{m}$ and $Z\mathbf{e}_{0}\left(
y,\xi\right)  $ vanish at the sharp point. Moreover, we have $\left|
Z\mathbf{e}_{0}\left(  y,\xi\right)  \right|  \leq c\left|  z-1/2-\imath
/2\right|  .$ Therefore the function $g$ can not be square integrable unless
$b_{0}=0.$ Similarly, the inclusion $Ag\in L_{2}\left(  \mathbb{Q}\right)  $
implies that $b_{1}=0$ and so on. Therefore $b_{0}=...=b_{m}=0$ and $g=0;$
which yields $c_{\lambda}=0$ for all $\lambda\in\Lambda.$ $\blacktriangleright$

\section{Gabor transform of a Gabor series}

Let $G$ be a set in the phase plane $\Phi$ and $r>0;$ for any $r>0$ we denote
by $G\left(  r\right)  $ the $r$-neighborhood of $G.$ Consider a convergent
series%
\[
g=\sum_{\lambda\in G}c^{\lambda}\mathbf{e}_{\lambda}%
\]
and estimate the Gabor transform of $g:$

\begin{proposition}
\label{t}For any subset $G\subset\Lambda$ and any $r>0$ we have%
\[
\int_{\Phi\backslash G\left(  r\right)  }\left|  \left\langle g|\mathbf{e}%
_{\mu}\right\rangle \right|  ^{2}\mathrm{d}\mu\leq\exp\left(  -\pi
r^{2}\right)  \sum\left|  c^{\lambda}\right|  ^{2}%
\]
\end{proposition}

\textsc{Proof. }By (\ref{ee}) we have for $\mu\in\Phi\backslash G\left(
r\right)  $%
\[
\left|  \left\langle g|\mathbf{e}_{\mu}\right\rangle \right|  ^{2}\leq\left(
\sum_{\lambda\in G}\left|  \left\langle \mathbf{e}_{\mu}|\mathbf{e}_{\lambda
}\right\rangle \right|  \left|  c^{\lambda}\right|  \right)  ^{2}\leq\sum
_{G}\exp\left(  -\pi\left|  \lambda-\mu\right|  ^{2}\right)  \sum_{G}\left|
c^{\lambda}\right|  ^{2}.
\]
Integrating the right-hand side yields
\[
\int_{\Phi\backslash G\left(  r\right)  }\left|  \left\langle g|\mathbf{e}%
_{\mu}\right\rangle \right|  ^{2}\mathrm{d}\mu\leq\sum_{\lambda\in G}\left|
c^{\lambda}\right|  ^{2}\int_{\Phi}\exp\left(  -\pi\left|  \lambda-\mu\right|
^{2}\right)  \phi\left(  \mu\right)  \mathrm{d}\mu.
\]
where $\phi$ is the indicator function of the set $\Phi\backslash G\left(
r\right)  .$ Write the sum over $\lambda$ as the integral over $\Phi$ with the
singular measure $\sigma\left(  \lambda\right)  =\sum_{G}\left|  c^{\lambda
}\right|  ^{2}\delta_{\lambda}.$ and change the variable $\lambda$ by
$\kappa=\lambda-\mu$ in the double integral%
\begin{align*}
\sum_{\lambda\in G}\int_{\Phi}\exp\left(  -\pi\left|  \lambda-\mu\right|
^{2}\right)  h\left(  \mu\right)  \mathrm{d}\mu\left|  c^{\lambda}\right|
^{2}  &  =\int_{\Phi}\int_{\Phi}\exp\left(  -\pi\left|  \lambda-\mu\right|
^{2}\right)  \phi\left(  \mu\right)  \mathrm{d}\mu\sigma\left(  \lambda\right)
\\
&  =\int_{\Phi}\exp\left(  -\pi\left|  \kappa\right|  ^{2}\right)  \int_{\Phi
}\phi\left(  \mu\right)  \sigma\left(  \kappa+\mu\right)  \mathrm{d}\mu
\end{align*}
The interior integral in the right-hand side vanishes if $\left|
\kappa\right|  <r,$ otherwise it is bounded by the total integral%
\[
\int\sigma\left(  \kappa+\mu\right)  \mathrm{d}\mu=\sum\left|  c^{\lambda
}\right|  ^{2}\mathrm{d}\kappa
\]
This yields the estimate for the right-hand side%
\[
\sum\left|  c^{\lambda}\right|  ^{2}\int_{\left|  \kappa\right|  \geq r}%
\exp\left(  -\pi\left|  \kappa\right|  ^{2}\right)  \mathrm{d}\kappa
=\sum\left|  c^{\lambda}\right|  ^{2}\exp\left(  -\pi r^{2}\right)
.\;\blacktriangleright
\]

\section{A `certainty' theorem}

\begin{theorem}
Let $K\subset D$ be arbitrary bounded domains in $\Phi$, such that $D$ is
$r$-neighborhood of $K$ for some $r\geq r_{0}$. An arbitrary function $f\in H$
can be written in the form
\[
f=\sum_{\lambda\in\Lambda\cap D}\alpha^{\lambda}\mathbf{e}_{\lambda}+\sum
_{\mu\in\left(  \Lambda+\sharp\right)  \cap D\backslash K}\omega^{\mu
}\mathbf{e}_{\mu}+\phi_{r}%
\]
where for $f\in H^{\delta}$%
\begin{equation}
\left\|  \phi_{r}\right\|  \leq\left(  \int_{\Phi\backslash D}\left|
\left\langle f|\mathbf{e}_{\mu}\right\rangle \right|  ^{2}\mathrm{d}%
\mu\right)  ^{1/2}+C_{\delta}r^{\delta}\exp\left(  -r/e\right)  \left\|
f\right\|  _{\delta}. \label{fi}%
\end{equation}
and the constant $C_{\delta}$ does not depend on $K$ and $D$.
\end{theorem}

\textbf{Remark 1.} Note that
\[
\alpha^{\lambda}=\gamma^{\lambda}\left(  f\right)  +O\left(  r^{\delta}%
\exp\left(  -r/e\right)  \right)
\]
for $\lambda\in K.$ The geometry of the two first terms in shown in Fig.1. \medskip

\begin{center}
\begin{picture}(400,400)
\multiput(0,40)(5,0){76}{\multiput(0,0)(0,5){68}{\circle*{2}}}
\put(380,200){\vector(1,0){20}} \put(385,205){$x$}
\put(200,380){\vector(0,1){20}}\put(190,390){$\xi$}
\qbezier(180,82)(-110,140)(170,280)
\qbezier(180,82)(430,60)(270,168)
\qbezier(240,240)(216,200)(270,168)
\qbezier(240,240)(330,350)(170,280) \thicklines
\qbezier(180,72)(-146,130)(170,294)
\qbezier(180,72)(450,46)(302,158)
\qbezier(260,242)(226,205)(302,158)
\qbezier(260,242)(350,370)(170,294) \thinlines\put(0,20){Fig.1:
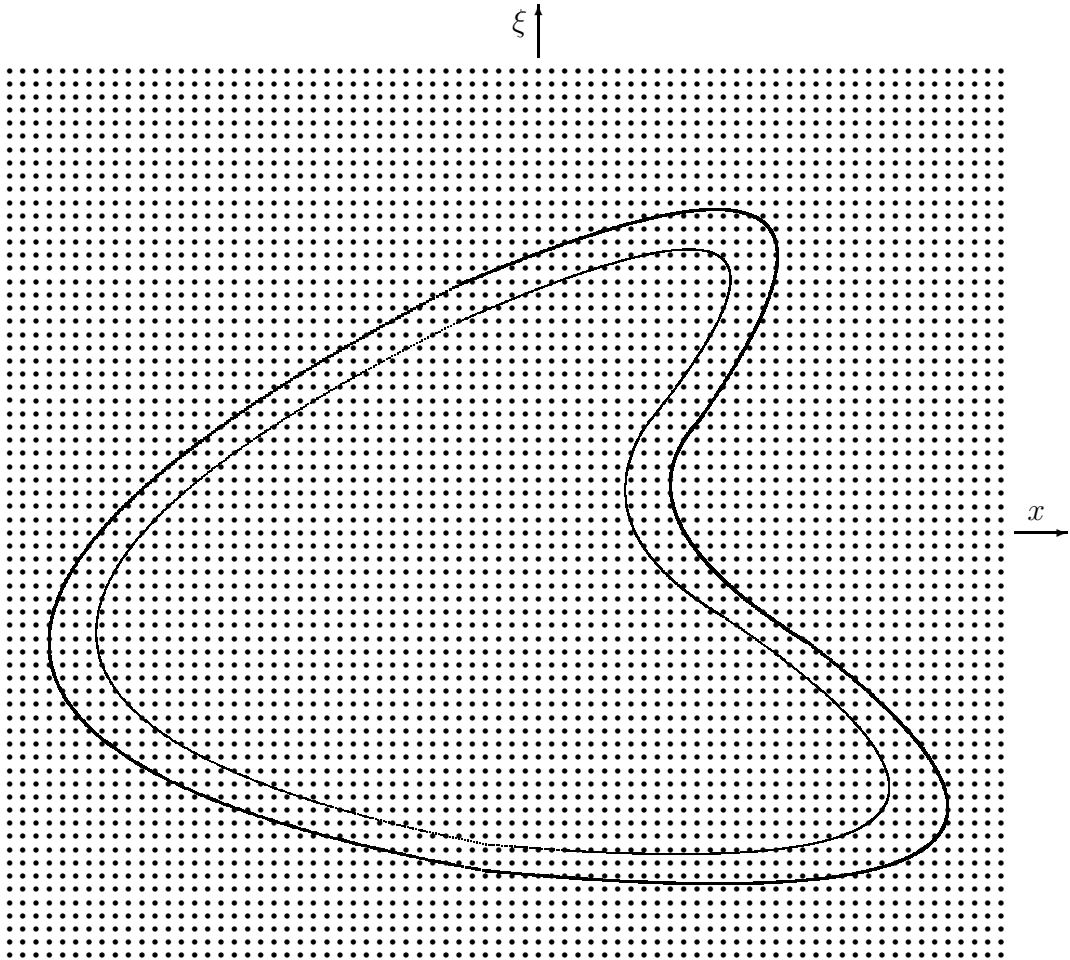
The lattice $\Lambda$, domains $K$ (thin lines) and $D$ (thick
lines) in the phase plane.}
\end{picture}
\end{center}

\textbf{Remark 2.} The exponential term in the estimate (\ref{fi}) is
indispensable, but its form might be made sharper. Indeed, vanishing of Gabor
means $\left\langle f|\mathbf{e}_{\mu}\right\rangle $ for $\mu\in
\Phi\backslash D$ does not guarantee that $f$ can be represented by a Gabor
functions supported in $\Lambda\cap D$ without additional term. This follows
from Proposition \ref{t} where the exponential factor $\exp\left(  -\pi
r^{2}\right)  $ is sharp.

\textsc{Proof of Theorem.} Set $K_{+}=K\left(  l\right)  ,$ $U=K\left(
r/2\right)  ,\;D_{-}=K\left(  r-l\right)  $ where the parameter $l=O\left(
r^{1/2}\right)  $ will be specified later. We have $K\subset K_{+}\subset
U\subset D_{-}\subset D.$ By Theorem \ref{t1} we can write $f=f_{U}+g,$ where%
\begin{equation}
f_{U}\doteq\sum_{\Lambda^{\sharp}\cap U}\gamma^{\lambda}\left(  f\right)
\mathbf{e}_{\lambda},\text{\ }g\doteq\sum_{\Lambda\backslash U}\gamma
^{\lambda}\left(  f\right)  \mathbf{e}_{\lambda}. \label{fg}%
\end{equation}
and the sharp point $\sharp$ is chosen in $U\backslash K.$ Estimate Gabor
transform of $g.$ By Proposition \ref{t} and Theorem \ref{t1}
\[
\int_{K_{+}}\left|  \left\langle g|\mathbf{e}_{\mu}\right\rangle \right|
^{2}\mathrm{d}\mu\leq\exp\left(  -\pi\left(  r/2-l\right)  ^{2}\right)
\sum_{\lambda\in\Lambda\backslash U}\left|  \gamma^{\lambda}\left(  f\right)
\right|  ^{2}\leq C\exp\left(  -\pi\left(  r/2-l\right)  ^{2}\right)  \left\|
f\right\|  _{\delta}^{2}%
\]
If $\mu\in\Phi\backslash D_{-}$, we can write $\left\langle g|\mathbf{e}_{\mu
}\right\rangle =\left\langle f|\mathbf{e}_{\mu}\right\rangle -\left\langle
f_{U}|\mathbf{e}_{\mu}\right\rangle $ since of (\ref{fg}). Apply again
Proposition \ref{t} and Theorem \ref{t1} and obtain
\[
\int_{\Phi\backslash D_{-}}\left|  \left\langle g|\mathbf{e}_{\mu
}\right\rangle \right|  ^{2}\mathrm{d}\mu\leq\int_{\Phi\backslash U_{+}%
}\left|  \left\langle f|\mathbf{e}_{\mu}\right\rangle \right|  ^{2}%
\mathrm{d}\mu+\exp\left(  -\pi\left(  r/2-l\right)  ^{2}\right)  \left\|
f\right\|  _{\delta}^{2},
\]
since $f_{U}$ is a sum of functions $\mathbf{e}_{\lambda},\lambda\in U.$ By
Proposition \ref{fee} we can write%
\begin{align}
g  &  =\int_{K_{+}}\left\langle g|\mathbf{e}_{\mu}\right\rangle \mathbf{e}%
_{\mu}\mathrm{d}\mu+\int_{\Phi\backslash D_{-}}\left\langle g|\mathbf{e}_{\mu
}\right\rangle \mathbf{e}_{\mu}\mathrm{d}\mu+\int_{D_{-}\backslash K_{+}%
}\left\langle g|\mathbf{e}_{\mu}\right\rangle \mathbf{e}_{\mu}\mathrm{d}%
\mu\label{ggg}\\
&  \doteq g_{+}+g_{-}+g_{0}\nonumber
\end{align}
The norms of the first two terms are bounded as follows
\begin{align}
\left\|  g_{+}\right\|  ^{2}  &  \leq\exp\left(  -\pi\left(  r/2-l\right)
^{2}\right)  \left\|  f\right\|  _{\delta}^{2},\;\label{gm}\\
\left\|  g_{-}\right\|  ^{2}  &  \leq\exp\left(  -\pi\left(  r/2-l\right)
^{2}\right)  \left\|  f\right\|  _{\delta}^{2}+\int_{\Phi\backslash D}\left|
\left\langle f|\mathbf{e}_{\mu}\right\rangle \right|  ^{2}\mathrm{d}%
\mu\label{gp}%
\end{align}

We fix a natural $m\leq r-1$ and transform the integral $g_{0}$ as follows:
for a point $\mu\in D_{-}\backslash K_{+}$ we choose a point $\lambda
\in\Lambda\cap D\backslash K$ and such that $\left|  \lambda-\mu\right|
\leq1/\sqrt{2}$ and consider the closed square $Q\left(  \mu\right)  $
centered at $\lambda$ with side $\sqrt{m+1}$, see Fig.2.

\begin{center}
\begin{picture}(300,240)
\multiput(60,40)(30,0){6}{\multiput(0,0)(0,30){6}{\multiput(8,9)(0,2){2}
{\line(1,0){4}} \multiput(9,8)(2,0){2}{\line(0,1){4}} }}
\put(60,40){\line(1,0){170}}\put(60,40){\line(0,1){170}}
\put(60,210){\line(1,0){170}}\put(230,40){\line(0,1){170}}
\put(152,130){\circle{3}} \put(148,135){$\mu$}
\put(145,125){\circle*{4}}\put(136,125){$\lambda$}
\multiput(85,65)(30,0){5}{\multiput(0,0)(0,30){5}{\circle*{3}}}
\put(0,0){Points of the set $(\Lambda+\sharp)\cap Q(\mu)$ are
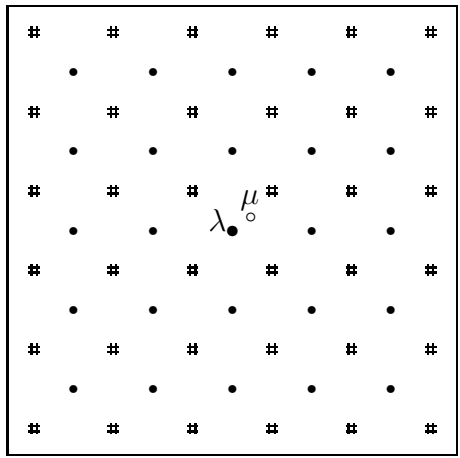
shown by $\sharp$;} \put(0,15){Fig.2: $\lambda$ is the closest
point to $\mu$ in the lattice $\Lambda$ (thick points) . }
\end{picture}
\end{center}

The set $Q\left(  \mu\right)  $ contains $m+1$ different points $\mu
_{0},...,\mu_{m}\in\Lambda+\sharp$ and is contained in the $l$-neighborhood of
$\mu$ where $l=\left(  \left(  m+1\right)  /2\right)  ^{1/2}+1.$ It follows
that $Q\left(  \mu\right)  \subset D\backslash K$ and
\begin{equation}
\max\left|  \mu-\mu_{s}\right|  \leq l,\;\max\left|  \lambda-\mu_{s}\right|
\leq l. \label{ll}%
\end{equation}
Applying Theorem \ref{m} to the function $\phi=\mathbf{e}_{\mu-\lambda}$ and
sharp points $\mu_{0}-\lambda,...,\mu_{m}-\lambda$ yields%
\begin{align}
\mathbf{e}_{\mu-\lambda}\left(  x\right)   &  =\sum_{j=0}^{m}\gamma^{\sharp
}\left(  \mathbf{a}^{k}\mathbf{e}_{\mu-\lambda}\right)  \mathbf{d}%
_{m,k}+g_{\mu,\lambda},\label{eee}\\
g_{\mu,\lambda}  &  \doteq\sum_{\Lambda}\gamma_{m}^{\nu}\left(  \mathbf{e}%
_{\mu-\lambda}\right)  \mathbf{e}_{\nu},\nonumber
\end{align}
where $\mathbf{d}_{j}$ belong to the linear span of $\mathbf{e}_{\mu
_{s}-\lambda},j=0,...,m.$ We have%
\begin{equation}
\left\|  g_{\mu,\lambda}\right\|  _{\delta,m}\leq\left\|  \mathbf{e}%
_{\mu-\lambda}\right\|  _{\delta,m}+\left\|  \sum\gamma^{\sharp}\left(
\mathbf{a}^{k}\mathbf{e}_{\mu-\lambda}\right)  \mathbf{d}_{m,k}\right\|
_{\delta,m}.\nonumber
\end{equation}
and $\mathbf{a}^{k}\mathbf{e}_{\mu-\lambda}=\left(  \mathbf{\mu}%
-\mathbf{\lambda}\right)  ^{k}\mathbf{e}_{\mu-\lambda}.$ By Proposition
\ref{d} and (\ref{ll})%
\begin{align*}
\left\|  \sum_{k}\gamma^{\sharp}\left(  \mathbf{a}^{k}\mathbf{e}_{\mu-\lambda
}\right)  \mathbf{d}_{m,k}\right\|  _{\delta,m}  &  \leq C\left(  l+1\right)
^{m+\delta}l^{m}\leq C\left(  m+1\right)  ^{m+\delta}\\
\left\|  \mathbf{e}_{\mu-\lambda}\right\|  _{\delta,m}  &  \leq C\left(
l+1\right)  ^{m+2}\leq C^{\prime}\left(  m+1\right)  ^{\left(  m+\delta
\right)  /2}%
\end{align*}
where the constants does not depend on $m$. This yields%
\begin{align*}
\left\|  \sum\gamma^{\sharp}\left(  \mathbf{a}^{k}\mathbf{e}_{\mu-\lambda
}\right)  \mathbf{d}_{m,k}\right\|  _{\delta,m}  &  \leq C\left(  m+1\right)
^{m+\delta},\\
\left\|  g_{\mu,\lambda}\right\|  _{\delta,m}  &  \leq C^{\prime}\left(
m+1\right)  ^{m+\delta}%
\end{align*}
By Theorem \ref{m}%
\begin{equation}
\sum\left(  \left|  \nu\right|  ^{2}+1\right)  ^{m}\left|  \gamma_{m}^{\nu
}\left(  \mathbf{e}_{\mu-\lambda}\right)  \right|  ^{2}\leq C\left(
m+1\right)  ^{2\left(  m+\delta\right)  }\left\|  \mathbf{e}_{\mu-\lambda
}\right\|  _{\delta,m}^{2}, \label{gmm}%
\end{equation}
Let $\lambda=\left(  p,\theta\right)  \in\Lambda;$ apply the operator
$T_{\lambda}$ as in (\ref{WH}). For $\mu=\left(  q,\xi\right)  ,\nu=\left(
r,\eta\right)  $ we obtain%
\[
T_{\lambda}\mathbf{e}_{\mu-\lambda}=\exp\left(  -2\pi\imath p\left(
\xi-\theta\right)  \right)  \mathbf{e}_{\mu},\;T_{\lambda}\mathbf{e}_{\nu
}=\exp\left(  -2\pi\imath p\eta\right)  \mathbf{e}_{\nu+\lambda}%
\]
and%
\begin{align}
\mathbf{e}_{\mu}\left(  x\right)   &  =\sum_{j=0}^{m}\gamma^{\sharp}\left(
\mathbf{a}^{j}\mathbf{e}_{\mu-\lambda}\right)  \exp\left(  2\pi\imath p\left(
\xi-\theta\right)  \right)  \mathbf{d}_{j,m}\label{em}\\
&  +\sum_{\Lambda}\gamma_{m}^{\nu}\left(  \mathbf{e}_{\mu-\lambda}\right)
\exp\left(  2\pi\imath p\left(  \xi-\theta+\eta\right)  \right)
\mathbf{e}_{\nu+\lambda}\nonumber
\end{align}
The function $\mathbf{d}_{j,m}$ belongs to the linear envelope of
$\mathbf{e}_{\kappa},\kappa\in\Lambda+\sharp$. Change the variable
$\nu+\lambda$ by $\nu$ in the second sum and write this equation in the form%
\begin{equation}
\mathbf{e}_{\mu}\left(  x\right)  =\sum_{\kappa\in\left(  \Lambda
+\sharp\right)  \cap Q\left(  \mu\right)  }\beta_{\mu}^{\kappa}\mathbf{e}%
_{\kappa}+\sum_{\lambda\in\Lambda\cap D}\varepsilon_{\mu}^{\lambda}%
\mathbf{e}_{\lambda}+\sum_{\nu\in\Lambda\backslash D}\omega_{\mu}^{\nu
}\mathbf{e}_{\nu} \label{ems}%
\end{equation}
where $\omega_{\mu}^{\nu}=\gamma_{m}^{\nu-\lambda}\left(  \mathbf{e}%
_{\mu-\lambda}\right)  \exp\left(  2\pi\imath r\left(  \xi-\theta+\eta\right)
\right)  .$ Estimate the third term by means of (\ref{gmm})%
\begin{align*}
\sum_{\nu}\left|  \omega_{\mu}^{\nu}\right|  ^{2}  &  =\sum_{\lambda\in
\Lambda\backslash D}\left|  \gamma_{m}^{\nu-\lambda}\left(  \mathbf{e}%
_{\mu-\lambda}\right)  \exp\left(  2\pi\imath r\left(  \xi-\theta+\eta\right)
\right)  \right|  ^{2}\\
&  \leq C\left(  m+1\right)  ^{2m+2}\left(  \left|  \nu-\lambda\right|
+1\right)  ^{-2m}\leq C\left(  m+1\right)  ^{2\left(  m+\delta\right)
}r^{-2m}%
\end{align*}
since $\left|  \nu-\lambda\right|  \geq r-2^{1/2}.$ Taking $m=r/e-1$ we obtain
$\left(  m+1\right)  ^{m+\delta}r^{-m}\leq Cr^{\delta}\exp\left(  -r/e\right)
.$ The relation $l=O\left(  m^{1/2}\right)  =O\left(  r^{1/2}\right)  $ is
fulfilled for this choice of $m.$ Taking in account that the square $Q\left(
\mu\right)  $ is contained in $D\backslash K$ $\ $and integrating (\ref{ems})
on $D_{-}\backslash K_{+}$ against the density $\left\langle g|\mathbf{e}%
_{\mu}\right\rangle \mathrm{d}\mu,$ yields%
\[
g_{0}=\sum_{\kappa\in\left(  \Lambda+\sharp\right)  \cap D\backslash K}%
\beta^{\kappa}\mathbf{e}_{\kappa}+\sum_{\lambda\in\Lambda\cap D}%
\varepsilon^{\lambda}\mathbf{e}_{\lambda}+\sum_{\lambda\in\Lambda\backslash
D}\omega^{\nu}\mathbf{e}_{\nu},
\]
where
\[
\beta^{\kappa}=\int_{D_{-}\backslash K_{+}}\beta_{\mu}^{\kappa}\left\langle
g|\mathbf{e}_{\mu}\right\rangle \mathrm{d}\mu,\;\varepsilon^{\lambda
}=...,\,\omega^{\nu}=...
\]
By Proposition \ref{p1} and Theorem \ref{t1} we have
\begin{align}
\left(  \sum\left|  \omega^{\nu}\right|  ^{2}\right)  ^{1/2}  &  \leq
C\exp\left(  -r/e\right)  \left(  \int\left|  \left\langle g|\mathbf{e}_{\mu
}\right\rangle \right|  ^{2}\mathrm{d}\mu\right)  ^{1/2}\leq Cr^{\delta}%
\exp\left(  -r/e\right)  \left\|  g\right\| \\
&  \leq C_{\delta}r^{\delta}\exp\left(  -r/e\right)  \left\|  f\right\|
_{\delta} \label{do}%
\end{align}
Finally we get%
\[
f=\sum_{\lambda\in\Lambda\cap D}\alpha^{\lambda}\mathbf{e}_{\lambda}%
+g_{+}+g_{-}+\sum_{\nu\in\left(  \Lambda+\sharp\right)  \cap D\backslash
K}\beta^{\nu}\mathbf{e}_{\nu}+\sum_{\lambda\in\Lambda\backslash D}\omega^{\nu
}\mathbf{e}_{\nu}%
\]
where $\alpha^{\lambda}=\gamma^{\lambda}\left(  f\right)  +\varepsilon
^{\lambda}$ for $\lambda\in\Lambda\cap U$ and $\alpha^{\lambda}=\varepsilon
^{\lambda}$ for $\lambda\in\Lambda\cap D\backslash U$ and the term
$\gamma^{\sharp}\left(  f\right)  \mathbf{e}_{\sharp}$ is included in the
second sum. We arrange this sum as follows%
\[
f=\sum_{\Lambda\cap D}\alpha^{\lambda}\mathbf{e}_{\lambda}+\sum_{\kappa
\in\left(  \Lambda+\sharp\right)  \cap D\backslash K}\beta^{\kappa}%
\mathbf{e}_{\kappa}+\phi_{r}%
\]
where and
\[
\phi_{r}=g_{+}+g_{-}+\sum_{\lambda\in\Lambda\backslash D}\omega^{\lambda
}\mathbf{e}_{\lambda}.\;\blacktriangleright
\]
By (\ref{gm}), (\ref{gp}) and (\ref{do}) we have%
\begin{align*}
\left\|  \phi_{r}\right\|   &  \leq\left\|  g_{+}\right\|  +\left\|
g_{-}\right\|  +\left\|  \sum\omega^{\nu}\mathbf{e}_{\nu}\right\| \\
&  \leq\left(  \int_{\Phi\backslash D}\left|  \left\langle f|\mathbf{e}_{\mu
}\right\rangle \right|  ^{2}\mathrm{d}\mu\right)  ^{1/2}+2\exp\left(
-\pi\left(  r/2-l\right)  ^{2}/2\right)  \left\|  f\right\|  _{\delta
}+C_{\delta}r^{\delta}\exp\left(  -r/e\right)  \left\|  f\right\|  _{\delta},
\end{align*}
which yields (\ref{fi}) for any sufficiently large $r$. $\blacktriangleright$

\section{Metaplectic group}

Remind that the Weyl-Heisenberg group is the space $X\times X^{\ast}%
\times\mathbb{R}$ with the group operation
\[
\left(  x,\xi,\tau\right)  \cdot\left(  x^{\prime},\xi^{\prime},\tau^{\prime
}\right)  =\left(  x+x^{\prime},\xi+\xi^{\prime},\tau+\tau^{\prime}+\frac
{1}{2}\left(  x^{\prime}\xi-x\xi^{\prime}\right)  \right)  .
\]
\textbf{Definition.} A linear transform $S$ of $\Phi$ is called
\textit{simplectic, }if it preserves the canonical bilinear form
$\sigma\left[  \left(  x,\xi\right)  ,\left(  y,\eta\right)  \right]  =\eta
x-\xi y.$ We can see that $\sigma\left[  u,v\right]  =\left\langle
u|Jv\right\rangle ,$ where
\[
J=\left(
\begin{array}
[c]{cc}%
0 & 1\\
-1 & 0
\end{array}
\right)
\]
In particular, $\lambda=\left(  p,\theta\right)  \mapsto J\lambda=\left(
\theta,-p\right)  $ is a linear simplectic transform. In the case
$X=\mathbb{R}$ a linear transformation $S$ in $\mathbb{R\times R}^{\ast}$ is
simplectic if the matrix $\left(  a,b,c,d\right)  $ of $S$ satisfies
$ad-bc=1$. For an arbitrary linear simplectic transformation $S$ such that
$b\neq0$ the integral transform
\[
M_{S}f\left(  x\right)  =\left(  \imath b\right)  ^{-1/2}\int\exp\left(
\pi\imath\left(  \frac{d}{b}x^{2}-\frac{2}{b}yx+\frac{a}{b}y^{2}\right)
\right)  f\left(  y\right)  \mathrm{d}y.
\]
is well defined operator in $L_{2}$. If $n=0,$ the operator
$M_{S}$ is defined as composition $M_{T}F$, where $F$ is the
Fourier transform and $T=\imath ^{1/2}SJ^{-1}$. The operator
$M_{S}$ has unitary closure in $L_{2}$. The equation
$M_{S}M_{T}=\pm M_{ST}$ holds for any simplectic transformations
$S,T$. It is called metaplectic (two-valued) representation of
$S.$ In particular, the metaplectic operator $M_{J}$ is the
Fourier transform $F$ times the factor $\imath^{1/2}.$ This
representation is single-valued on the irreducible two-fold
covering of the group $S.$ See more information in \cite{Fo}.

Any rotation $S(x,\xi)=(\cos\varphi\,x-\sin\varphi\,\xi,\sin\varphi
\,x+\cos\varphi\,\xi)$ is an orthogonal simplectic transformation.

\begin{proposition}
If $S$ is a rotation as above, then
\begin{align*}
M_{S}\,\mathbf{a}  &  =\exp\left(  -\imath\varphi\right)  \,\mathbf{a\,}%
M_{S}\\
M_{S}\,\mathbf{a}^{+}  &  =\exp\left(  \imath\varphi\right)  \,\mathbf{a}%
^{+}M_{S}%
\end{align*}
where $\mathbf{a}^{+},$ $\mathbf{a}$ is the creation and the annihilation
operator, respectively.
\end{proposition}

\textsc{Proof.} Direct calculation.

\begin{proposition}
\label{rot}For any rotation $S$ in the phase plane and any point $\lambda$\ we
have
\begin{equation}
M_{S}\left(  \mathbf{e}_{\lambda}\right)  =\pm\exp\left(  \imath
\varphi/2\right)  \exp\left(  \pi\imath\left(  p\theta-q\eta\right)  \right)
\mathbf{e}_{S\left(  \lambda\right)  } \label{me}%
\end{equation}
where $\left(  q,\eta\right)  =S\left(  p,\theta\right)  $
\end{proposition}

\textsc{Proof }is straightforward.

\textbf{Remark. }The equation (\ref{me}) means that $M_{S}$ transforms a Gabor
function $\mathbf{e}_{\lambda}$ to another Gabor function (up to a phase
factor) while the `quantum support' $\lambda$ of a Gabor function moves by
action of the corresponding geometric transform $S.$ In particular, the
Fourier transform belongs to the metaplectic group: $F=M_{J}.$ We have for any
$\lambda=\left(  p,\theta\right)  ,$ $\mathbf{\hat{e}}_{\lambda}(\eta)=\pm
\exp\left(  \imath\varphi/2\right)  \mathbf{e}_{\widehat{\lambda}}\left(
\eta\right)  ,\,\widehat{\lambda}=J\lambda=\left(  \theta,-p\right)  $ for
a\ real $\phi.$ The metaplectic representation can be thought as a
`quantization' of group of simplectic transforms.

\begin{corollary}
For any $\delta>1$ and natural $m$ the space $H^{\delta,m}$ is invariant under
action of the metaplectic representation of the rotation group.
\end{corollary}

\textsc{Proof.} For any rotation $S$ we have
\[
\left\langle f|\mathbf{e}_{\lambda}\right\rangle =\left\langle M_{S}%
f|M_{S}\mathbf{e}_{\lambda}\right\rangle =\exp\left(  \imath\varphi/2\right)
\left\langle M_{S}f|\mathbf{e}_{\lambda}\right\rangle
\]
which yields $\left|  \left\langle M_{S}f|\mathbf{e}_{\lambda}\right\rangle
\right|  =\left|  \left\langle f|\mathbf{e}_{\lambda}\right\rangle \right|  .$
It follows that $H^{\delta}$ is invariant. The same true for $H^{\delta,m}$
since of (\ref{me}).\ $\blacktriangleright$

\begin{corollary}
For an arbitrary orthogonal simplectic transformation $S$ in the phase space,
theorem \ref{t1} holds for functions $f\in H$ and Gabor system $\mathbf{e}%
_{\lambda},\lambda\in$ $S\left(  \Lambda^{\sharp}\right)  .$
\end{corollary}

\textbf{Remark}. The operator $M$ relates to the metaplectic
representation $\mu$ in the sense of \cite{Fo} by the equation
$M_{S}=\mu\left(  \tilde {S}\right)  ,$ where the simplectic
matrix $\tilde{S}$ is obtained from $S$ by changing sign at $b$
and $c$ and replacing $\theta$ by $-\theta.$ This results in the
representation $\mu$ which is chosen in \cite{Fo}. Note that the
mapping $S\mapsto\tilde{S}$ is a involution in the simplectic
group.

\end{document}